\documentclass[12pt]{article}
\usepackage{amssymb,amsfonts}
\usepackage{amsmath,amsthm,amsxtra}
\numberwithin{equation}{section}
\usepackage{bbm}
\usepackage[all]{xy}
\usepackage{graphicx}
\usepackage{color}
\usepackage{verbatim}
\usepackage[notcite,notref]{showkeys}

\usepackage{enumerate}

\usepackage[utf8]{inputenc}  
\usepackage[T1]{fontenc}

\newcommand\bigcheck[1]{#1 \raise1ex\hbox{$\hspace{-1ex}{}^\vee$}}
\newcommand\sucheck[1]{#1 \raise0.5ex\hbox{$\hspace{-1ex}{}^\vee$}}

\makeatletter

\renewcommand\section{\@startsection {section}{ 0}{\z@}%
                                    {-4ex \@plus -2ex \@minus -.2ex}%
                                    {6ex \@plus.2ex}%
                                    {\centering\large\scshape}}
\renewcommand\subsection{\@startsection{subsection}{1}{\z@}%
                                     {-5ex\@plus -6ex \@minus -3ex}%
                                     {-75ex \@plus 100ex}%
                                     {\normalfont\large\bfseries}}
\renewcommand\subsubsection{\@startsection{subsubsection}{3}{\z@}%
                {-3.25ex\@plus -1ex \@minus -.2ex}%
                {1.5ex \@plus .2ex}%
                {\normalfont\normalsize\bfseries}}
\makeatother

\makeatletter
\newcommand\xleftrightarrow[2][]{%
  \ext@arrow 9999{\longleftrightarrowfill@}{#1}{#2}}
\newcommand\longleftrightarrowfill@{%
  \arrowfill@\leftarrow\relbar\rightarrow}
\makeatother

\newtheorem{innercustomgeneric}{\customgenericname}
\providecommand{\customgenericname}{}
\newcommand{\newcustomtheorem}[2]{%
  \newenvironment{#1}[1]
  {%
   \renewcommand\customgenericname{#2}%
   \renewcommand\theinnercustomgeneric{##1}%
   \innercustomgeneric
  }
  {\endinnercustomgeneric}
}

\newcustomtheorem{customthm}{Theorem}
\newcustomtheorem{customlemma}{Lemma}
\newcustomtheorem{customremark}{Remark}

\newcommand{\mc}[1]{{\mathcal #1}}

\usepackage[pdftex]{hyperref}
\hypersetup{
    bookmarks=true,         
    unicode=false,          
    pdftoolbar=true,        
    pdfmenubar=true,        
    pdffitwindow=true,      
    pdftitle={My title},    
    pdfauthor={Author},     
    pdfsubject={Subject},   
    pdfnewwindow=true,      
    pdfkeywords={keywords}, 
    colorlinks=false,       
    linkcolor=red,          
    citecolor=green,        
    filecolor=magenta,      
    urlcolor=cyan           
}
\begin{document}
\begin{titlepage}
	\centering
	\vspace{3.5cm}
	{\Huge\bfseries Compatible Hamiltonian Operators for the Krichever-Novikov Equation\par}
	\vspace{2cm}
	{\Large\itshape Sylvain Carpentier*\par}
\vspace{1 cm}
\begin{abstract}
It has been proved by Sokolov that Krichever-Novikov equation's hierarchy is hamiltonian for the Hamiltonian
operator $H_0=u_x \partial^{-1} u_x$ and possesses two weakly non-local recursion operators of
 degree $4$ and $6$, $L_4$ and $L_6$. We show here that $H_0$, $L_4H_0$ and $L_6H_0$
are compatible Hamiltonians operators for which the Krichever-Novikov equation's hierarchy is hamiltonian. 
\end{abstract}
	{\large \today\par}
\vspace{1 cm}

{\let\thefootnote\relax\footnotetext{* \textit{Department of Mathematics, Massachusetts Institute of Technology, Cambridge, MA 02139, USA}}}
\end{titlepage}

In the study of finite gap solutions of KP, an integrable $(1+1)$-dimensional PDE was discovered, the Krichever-Novikov equation. One of its forms (equivalent to the original one in [KN80]) is
\begin{equation}\tag{1}
\frac{du}{dt}=u_3-\frac{3}{2} \frac{u_2^2}{u_1}+\frac{P(u)}{u_1},
\end{equation}
where $u=u(t,x)$, $u_n=(\frac{d}{dx})^n(u)$, and $P$ is a polynomial of degree at most $4$. Let $\mc V=\mathbb{C}[u,u_1^{\pm},u_2,...]$ and $\mc K$ be the fraction field of $\mc V$. Let us denote $\frac{d}{dx}$ by $\partial$. The \textit{differential order} $d_F$ of a function $F \in \mc V$ is the highest integer $n$ such that $\frac{\partial F}{\partial u_n} \neq 0$.
\par
 One of the attributes of equation $(1)$ is to be part of an infinite hierarchy of compatible evolution PDEs of odd differential orders 
\begin{equation}\tag{2}
\frac{du}{dt_i}=G_i \in \mc V, i \geq 0,
\end{equation}
where $G_i$ has differential order $(2i+1)$. One says that $F, G \in \mc V$ are \textit{compatible}, or \textit{symmetries} of one another, if 
\begin{equation}\tag{3}
\{F,G\}:=X_F(G)-X_G(F)=0,
\end{equation}
where $X_F$ denotes the derivation of $\mc V$ induced by the evolution equation $u_t=F$, that is
\begin{equation}\tag{4}
X_F=\sum_{n \geq 0}{F^{(n)}\frac{\partial}{\partial u_n}}.
\end{equation}
$(3)$ endows $\mc V$ with a Lie algebra bracket, and the $G_i$'s span an infinite-dimensional abelian subalgebra of $(\mc V, \{.,.\})$, which we will denote by $\mc S$.
The first four equations in the hierarchy are:
\begin{equation}\tag{5}
\begin{split}
G_0=& \hspace{1 mm}u_1,\\
G_1=&\hspace{1 mm}u_3-\frac{3}{2} \frac{u_2^2}{u_1}+\frac{P(u)}{u_1},\\
G_2= &\hspace{1 mm}u^{(5)}-5 \frac{u_4u_2}{u_1}-\frac{5}{2} \frac{u_3^2}{u_1}+ \frac{25}{2} \frac{u_3u_2^2}{u_1^2}-\frac{45}{8}\frac{u_2^4}{u_1^3}-\frac{5}{3}P \frac{u_3}{u_1^2}+\frac{25}{6} P \frac{u_2^2}{u_1^3}\\ &-\frac{5}{3}P_u\frac{u_2}{u_1}-
\frac{5}{18} \frac{P^2}{u_1^3}+\frac{5}{9}u_1P_{uu},
\end{split}
\end{equation}
\begin{equation}\tag{6}
\begin{split}
G_3=&\hspace{1 mm} u_7-7\frac{u_2u_6}{u_1}-\frac{7}{6}\frac{u_5}{u_1^2}(2P+12u_3u_1-27u_2^2)-\frac{21}{2}\frac{u_4^2}{u_1}+\frac{21}{2}\frac{u_4}{u_1^3}(2P-11u_2^2)\\
&-\frac{7}{3}\frac{u_4}{u_1^2}(2P_u u_1-51u_2u_3)+\frac{49}{2}\frac{u_3^3}{u_1^2}+\frac{7}{12}\frac{u_3^2}{u_1^3}(22P-417u_2^2)+\frac{2499}{8}\frac{u_2^4}{u_1^4}u_3 \\
&\frac{91}{3}P_u\frac{u_2}{u_1^2}u_3-\frac{595}{6}P\frac{u_2^2}{u_1^4}u_3-\frac{35}{18}\frac{u_3}{u_1^4}(2P_{uu}u_1^4-P^2)-\frac{1575}{16}\frac{u_2^6}{u_1^5}+\frac{1813}{24}\frac{u_2^4}{u_1^5}P\\
&-\frac{203}{6} \frac{u_2^3}{u_1^3}P_u+\frac{49}{36}\frac{u_2^2}{u_1^5}(6P_{uu} u_1^4-5P^2)-\frac{7}{9}\frac{u_2}{u_1^3}(2P_{uuu}u_1^4-5PP_u)+\frac{7}{54}\frac{P^3}{u_1^5}\\
& -\frac{7}{9}P_{uu}\frac{P}{u_1}+\frac{7}{9}P_{uuuu}u_1^3-\frac{7}{18}\frac{P_u^2}{u_1}.
\end{split}
\end{equation}

It is known ([IS80], [MS08]) that all integrable hierarchies admit a pseudodifferential operators $L \in \mc V((\partial^{-1}))$ satisfying
\begin{equation}\tag{7}
X_F(L)=[D_F,L]
\end{equation}
for all $F$ in the hierarchy, where $D_F$ denotes the \textit{Fréchet derivative} of $F$:
\begin{equation}\tag{8}
D_F=\sum_{n}{\frac{\partial F}{\partial u_n} \partial^n} \in \mc V[\partial].
\end{equation}
A pseudodifferential operator satisfying $(7)$ is called a \textit{recursion operator} (for $F$). In [DS08] two rational recursions operators for $(1)$ were found, of order $4$ and $6$:
\begin{equation}\tag{9}
L_4=H_1 H_0^{-1}, \hspace{3 mm} L_6=H_2 H_0^{-1},
\end{equation}
where
\begin{equation}\tag{10}
\begin{split}
H_0=& u_1 \partial^{-1} u_1, \\
H_1=& \frac{1}{2}(u_1^2 \partial^3+\partial^3 u_1^2)+(2u_3u_1-\frac{9}{2}u_2^2-\frac{2}{3}P)\partial+\partial (2u_3u_1-\frac{9}{2}u_2^2-\frac{2}{3}P)\\
&+G_1 \partial^{-1} G_1+u_1 \partial^{-1}G_2+G_2 \partial^{-1} u_1,\\
H_2=&\frac{1}{2}(u_1^2 \partial^5+\partial^5 u_1^2)+(3u_3u_1-\frac{19}{2}u_2^2-P)\partial^3+\partial^3 (3u_3u_1-\frac{19}{2}u_2^2-P)\\
&+(u_5u_1-9u_3u_2+\frac{19}{2}u_3^2-\frac{2}{3}\frac{u_3}{u_1}(5P-39u_2^2)+\frac{u_2^2}{u_1^2}(5P-9u_2^2)+\frac{2}{3}\frac{P^2}{u_1^2}+u_1^2P_{uu}) \partial \\
&+\partial (u_5u_1-9u_3u_2+\frac{19}{2}u_3^2-\frac{2}{3}\frac{u_3}{u_1}(5P-39u_2^2)+\frac{u_2^2}{u_1^2}(5P-9u_2^2)+\frac{2}{3}\frac{P^2}{u_1^2}+u_1^2P_{uu})\\
&+ G_1 \partial^{-1}G_2+G_2 \partial^{-1} G_1 +u_1 \partial^{-1} G_3+G_3 \partial^{-1} u_1.
\end{split}
\end{equation}
Moreover, $L_4$ and $L_6$ are both \textit{weakly non-local}, i.e. of the form
\begin{equation}\tag{11}
E(\partial) \in \mc V[\partial]+\sum_i{p_i \partial^{-1} \frac{\delta \rho_i}{\delta u}},
\end{equation}
where the $\rho_i$'s are conserved densities of $(1)$. Recall that the \textit{variational derivative} $\frac{\delta}{\delta u}$ is defined as follows:
\begin{equation}\tag{12}
\frac{\delta F}{\delta u}=D_F^*(1)=\sum_n{(-\partial)^n(\frac{\partial F}{\partial u_n})}.
\end{equation}
 In [S84], Sokolov showed that the space of symmetries of $(1)$, $\mc S$, is preserved by $L_4$. The same argument applies to
$L_6$, which was found later. He also establishes that the hierarchy of the Krichever-Novikov equation is \textit{hamiltonian} for $H_0$:
there exists a sequence $\phi_i \in \mc V$ such that
\begin{equation}\tag{13}
G_i=H_0(\frac{\delta \phi_i}{\delta u}) \text{  for all  } i \geq 0.
\end{equation}
A \textit{Hamiltonian operator} $H=AB^{-1} \in \mc V(\partial)$ with $A$ and $B$ right coprime is a skewadjoint rational differential operator inducing a non-local poisson lambda bracket, which is equivalent to the following identity (see equation $(6.13)$ in [DSK13])
\begin{equation}\tag{14}
\begin{split}
&A^*(D_{B(F)}(A(G))+D_{A(G)}^*(B(F))-D_{B(G)}(A(F))+D_{B(G)}^*(A(F)))\\
&=B^*(D_{A(G)}(A(F))-D_{A(F)}(A(G)) ) 
\end{split}
\end{equation}
for all $F,G \in \mc V$.
\begin{customlemma}{1}
Let $L \in \mc V(\partial)$ be a skewadjoint rational operator. If there exists an infinite-dimensional (over $\mathbb{C}$) subspace $\mc W \subset \mc V$ such that $B(W) \subset \frac{\delta}{\delta u} \mc V$ and such that for all $G \in \mc W$, $E=A(G)$ satisfies 
\begin{equation}\tag{15}
X_E(L)=D_EL+LD_E^*,
\end{equation}
then $L$ is a Hamiltonian operator. Conversely, if $L$ is a hamiltonian operator and $G \in \mc V$, then $D_{B(G)}=D_{B(G)}^*$ if and only if $A(G)$ satisfies equation $(15)$ 
\end{customlemma}
\begin{proof}
Let us first give an equivalent form of $(15)$ involving only differential operators. 
\begin{equation}\tag{16}
\begin{split}
(1.15) &\iff X_E(A)-D_EA=AB^{-1}(X_E(B)+D_E^*B)\\
&\iff X_E(A)-D_EA=-{B^*}^{-1}A^*(X_E(B)+D_E^*B)\\
&\iff A^*(X_E(B)+D_E^*B)=B^*(D_EA-X_E(A)) \\
&\iff A^*(X_E+D_E^*)B=B^*(D_E-X_E)A.\\
&\iff A^*(D_{B(F)}(E)+D_E^*(B(F)))=B^*(D_E(A(F))-D_{A(F)}(E)) \hspace{2 mm} \forall F \in \mc V
\end{split}
\end{equation}
Comparing the last line of $(16)$ with $(14)$, it is clear that if $H$ is Hamiltonian, then $E=A(G)$ satisfies equation $(15)$ is and only if $D_{B(G)}$ is self-adjoint. It is also clear that if $A(G)$ satisfies $(15)$ and $D_{B(G)}$ is self-adjoint, then $(F,G)$ satisfies $(14)$ for any $F \in \mc V$. Therefore, if we consider $\mc W \subset \mc V$ infinite-dimensional subspace of $\mc V$ such that $A(\mc W)$ satisfies $(15)$ and $B(\mc W) \subset \frac{\delta}{\delta u} \mc V$, we deduce that $(14)$ is satisfied for any $(F,G) \in \mc V \times \mc W$. To conclude, we note that $(14)$ can be rewritten as an identity of bidifferential operator, i.e. it amounts to say that some expression of the form $\sum m_{ij} F^{(i)}G^{(j)}$, where $m_{ij} \in \mc V$ is trivial, i.e. $m_{ij}=0$ for all $i,j$. Namely, $(14)$ is equivalent to
\begin{equation}\tag{17}
\begin{split}
&A^*(X_{A(G)}(B)(F)-X_{A(F)}(B)(G) +(D_A)_G^*(B(F))+(D_B)_G^*(A(F))) \\
&=B^*(X_{A(F)}(A)(G)-X_{A(G)}(A)(F)),
\end{split}
\end{equation}
where given a differential operator $P$, an element $F \in \mc V$, the differential operator $(D_P)_F$ is defined by
\begin{equation}\tag{18}
(D_P)_F(G)=X_G(P)(F) \hspace{2 mm} \forall G \in \mc V.
\end{equation}
If a bidifferential operator vanishes on $\mc V \times \mc W$, it must be identically $0$, since $\mc W$ is infinite dimensional. Hence, $L$ is an Hamiltonian operator.
\end{proof}

\begin{customlemma}{2}
Let $L=CD^{-1}$ be a rational operator and $(F_n)_{n \geq 0}$ a sequence spanning an infinite-dimensional subspace of $\mc K$ satisfying $C(F_n)=D(F_{n+1}) \in \mc V$ for all $n \geq 0$. Assume that $L$ is recursion for all the $D(F_n)$'s and that the $D(F_n)$'s are hamiltonian for some Hamiltonian operator $H \in \mc V(\partial)$. Then, provided $LH$ is skew-adjoint, $LH$ is a Hamiltonian operator for which all the $D(F_n)$'s are hamiltonian ($n \geq 1$).
\end{customlemma}
\begin{proof}
By Lemma $1$, $H$ satisfies equation $(15)$ for all $D(F_n), n \geq 0$, hence so does $LH$ ($L$ is recursion for $D(F_n)$ for all $n \geq 0$). To conclude using Lemma $1$, one needs to check that $D(F_n)=LH(\frac{\delta\rho_n}{\delta u})$ for some $\rho_n \in \mc V$ for all $n \geq 1$. Let $P,Q \in \mc V[\partial]$ be right coprime differential operators such that $LH=PQ^{-1}$. Let $A,B$ be right coprime differential operators such that $H=AB^{-1}$. $D(F_n)$ is hamiltonian for $H$ for all $n \geq 0$, meaning that there exist two sequences in $\mc V$, $(\phi_n)_{n \geq 0}$ and $(\rho_n)_{n \geq 0}$ , such that $\frac{\delta \rho_n}{\delta u}=B(\phi_n)$ and $D(F_n)=A(\phi_n)$ for all $n \geq 0$. In the language of [CDSK15], $\frac{\delta \rho_n}{\delta u}$ and $C(F_n)$ are $CD^{-1}AB^{-1}$ associated, hence (quote result) there exists $\psi_n$ such that $C(F_n)=P(\psi_n)$ and $Q(\psi_n)=\frac{\delta \rho_n}{\delta u}$ for all $n \geq 0$. Therefore, by Lemma $1.1$, $LH$ is a Hamiltonian operator for which $(C(F_n))_{n \geq 0}$ are hamiltonian.
\end{proof}

\begin{customthm}{}
$H_0$, $H_1$ and $H_2$ are compatible Hamiltonian operators.
\end{customthm}
\begin{proof}
Let $\alpha, \beta, \gamma \in \mathbb{C}$ and let $L_{\alpha,\beta,\gamma}=(\alpha H_0+\beta H_1+ \gamma H_2)H_0^{-1}$. $L_{\alpha,\beta,\gamma}$ is a recursion operator for the whole Krichever-Novikov hierarchy $\mc S$. Moreover, it maps $\mc S$ to itself as was proved in [S84], meaning that if $L_{\alpha,\beta,\gamma}=AB^{-1}$ with $A,B$ right coprime and $G \in \mc S$, then $G=B(F)$ for some $F \in \mc K$ and $A(F) \in \mc S$. The theorem follows from Lemma $2$.
\end{proof}

\begin{customremark}{3}
It follows from Lemma $1$ that $H=H_2H_1^{-1}H_0$ is a Hamiltonian operator of degree $1$. However, it is not weakly non-local. More generally all the $(H_2H_1^{-1})^nH_0$, for $n \in \mathbb{Z}$ are pairwise compatible Hamiltonian operators. 
\end{customremark}

\begin{customremark}{4}
Every Hamiltonian operator $K=AB^{-1}$ over $\mc V$, where $A$ and $B$ are right coprime induces a Lie algebra bracket on the space of functionals $\mc F(K) :=\{ \int f \in \mc V/ \partial \mc V| \frac{\delta f}{\delta u} \in Im B\}$, (well-)defined by $\{\int f, \int g\}=\int \frac{\delta f}{\delta u}AB^{-1}(\frac{\delta g}{\delta u})$ (see section $7.2$ in [DSK13]). Note that $\mc F(H_0)=\mc V/ \partial \mc V$ but that $\mc F(H_1)$ and $\mc F(H_2)$ consist only of the conserved densities of the Krichever-Novikov equation. 
\end{customremark}
We recall that if a rational differential operator $L=AB^{-1}$, with $A,B \in \mc V[\partial]$ right coprime generates an infinite dimensional abelian subspace of $(\mc V, \{.,.\})$, in the sense that there exist $(F_n)_{n \geq 0} \in \mc K$ such that $A(F_n)=B(F_{n+1})$ for all $n \geq 0$ and such that the $B(F_n)$'s span an infinite-dimensional abelian subspace of $(\mc V, \{.,.\})$, then for all $\lambda \in \mathbb{C}$, $Im (A + \lambda B)$ must be a sub Lie algebra of $(\mc V, \{.,.\})$ (see [C17]). The recursion operators $L_{\alpha,\beta,\gamma}$ satisfy this condition.
\par
The author thanks Vladimir Sokolov for useful discussions, and Victor Kac for his interest in this work.

\end{document}